\documentclass[12pt,leqno]{article}
\usepackage{amsmath, amsthm, amsfonts, amssymb, color}
\usepackage{mathrsfs}
\theoremstyle{definition}

\newcommand{\scr}[1]{\mathscr #1}
\definecolor{wco}{rgb}{0.5,0.2,0.3}

\numberwithin{equation}{section} \theoremstyle{remark}

\newcommand{\ua}{\uparrow}

\title{{\bf Derivative Formula and Applications  for Degenerate Diffusion Semigroups }\footnote{Supported in
 part by  NNSFC(11131003), SRFDP, 985 project through the Laboratory of Mathematical and  Complex Systems, and the Fundamental Research Funds for the Central Universities.}
}
\author{
{\bf Feng-Yu Wang$^{a),b)}$ and Xi-Cheng Zhang$^{c)}$ }\\
\footnotesize{a) School of Mathematical Sciences,
Beijing Normal University, Beijing 100875, China}\\
 \footnotesize{b) Department of Mathematics,
Swansea University, Singleton Park, SA2 8PP, UK}\\
\footnotesize{c) School of Mathematics and Statistics,
Wuhan University, Wuhan 430072, China}}
\begin{document}
\def\R{\mathbb R}  \def\ff{\frac} \def\ss{\sqrt} \def\B{\mathbf
B}
\def\N{\mathbb N} \def\kk{\kappa} \def\m{{\bf m}}
\def\dd{\delta} \def\DD{\Delta} \def\vv{\varepsilon} \def\rr{\rho}
\def\<{\langle} \def\>{\rangle} \def\GG{\Gamma} \def\gg{\gamma}
  \def\nn{\nabla} \def\pp{\partial} \def\EE{\scr E}
\def\d{\text{\rm{d}}} \def\bb{\beta} \def\aa{\alpha} \def\D{\scr D}
  \def\si{\sigma} \def\ess{\text{\rm{ess}}}
\def\beg{\begin} \def\beq{\begin{equation}}  \def\F{\scr F}
\def\Ric{\text{\rm{Ric}}} \def\Hess{\text{\rm{Hess}}}
\def\e{\text{\rm{e}}} \def\ua{\underline a} \def\OO{\Omega}  \def\oo{\omega}
 \def\tt{\tilde} \def\Ric{\text{\rm{Ric}}}
\def\cut{\text{\rm{cut}}} \def\P{\mathbb P} \def\ifn{I_n(f^{\bigotimes n})}
\def\C{\scr C}      \def\aaa{\mathbf{r}}     \def\r{r}
\def\gap{\text{\rm{gap}}} \def\prr{\pi_{{\bf m},\varrho}}  \def\r{\mathbf r}
\def\Z{\mathbb Z} \def\vrr{\varrho} \def\ll{\lambda}
\def\L{\scr L}\def\Tt{\tt} \def\TT{\tt}\def\II{\mathbb I}
\def\i{{\rm in}}\def\Sect{{\rm Sect}}\def\E{\mathbb E} \def\H{\mathbb H}
\def\M{\scr M}\def\Q{\mathbb Q} \def\texto{\text{o}} \def\LL{\Lambda}
\def\Rank{{\rm Rank}} \def\B{\scr B} \def\i{{\rm i}} \def\HR{\hat{\R}^d}

\maketitle
\begin{abstract} By using the Malliavin calculus and solving a control problem, Bismut type derivative formulae are established for a class of degenerate diffusion semigroups with non-linear drifts. As applications, explicit gradient estimates and Harnack inequalities are derived.
\end{abstract} \noindent

 \vskip 2cm

\section{Introduction}

The Bismut derivative formula  introduced in \cite{Bismut},  also known as   Bismut-Elworthy-Li formula  due to \cite{EL},  is a powerful tool to derive regularity estimates on diffusion semigroups. In the elliptic case this formula can be expressed by using the intrinsic curvature induced by the generator.
But in the degenerate case the required curvature lower bound is no longer available.
Of course, the Malliavin calculus   works also for the hypoelliptic case as shown in e.g. \cite{AT} on Riemannian manifolds.
In this case the pull-back operator involved in the formula is normally less explicit, so that it is hard for one to derive explicit gradient estimates. Nevertheless, as shown in \cite[\S 6]{AT}, in some concrete degenerate cases the derivative formula can be explicitly established by solving certain control problems.

Recently,  explicit derivative formulae  for  damping stochastic Hamiltonian systems have been established in \cite{Z} and \cite{GW} by using Malliavin calculus and coupling respectively, where the degenerate part is linear. In this case successful couplings with control can be constructed in a very explicit way, so that some  known arguments developed in the elliptic setting can be applied. However, when the degenerate part is non-linear, the study becomes much more complicated. The main purpose of this paper is to extend results derived in \cite{Z, GW} to the non-linear degenerate case.

Consider the following degenerate stochastic differential equation on $\R^m\times \R^d$:

\beq\label{E1} \beg{cases} \d X_t^{(1)}= Z^{(1)}(X_t^{(1)}, X_t^{(2)})\d t,\\
\d X_t^{(2)} = Z^{(2)}(X_t^{(1)}, X_t^{(2)})\d t +\si \d B_t,\end{cases}\end{equation} where
$X_t^{(1)}$ and $ X_t^{(2)}$ take values in $\R^m$ and $\R^d$ respectively, $\si$ is an invertible $d\times d$-matrix, $B_t$ is a $d$-dimensional Brownian motion,
$Z^{(1)} \in C^2(\R^{m+d}; \R^m)$ and  $ Z^{(2)}\in C^1(\R^{m+d}; \R^d).$   Let $X_t=(X_t^{(1)},X_t^{(2)}), Z= (Z^{(1)}, Z^{(2)}).$ Then the equation can be formulated as
\beq\label{E2}  \d X_t = Z(X_t)\d t + (0,\si \d B_t).\end{equation}
We assume that the solution is non-explosive, which is ensured by (H1) below.
Our purpose is to establish an explicit   derivative formula for the associated Markov semigroup $P_t$:
$$P_t f(x)= \E f(X_t(x)),\ \ t>0, x\in\R^{m+d}, f\in \B_b(\R^{m+d}),$$
where $X_t(x)$ is the solution of (\ref{E2}) with $X_0=x$, and $\B_b(\R^{m+d})$ is the set of all bounded measurable functions on $\R^{m+d}.$

When $m=d,\si=I_{d\times d}$ and
$$Z^{(1)}(x,y)=\nn H(x,\cdot)(y),\ \ Z^{(2)}(x,y)= -\nn H(\cdot,y)(x)-F(x,y)\nn H(x,\cdot)(y)$$ for some functions $H$ and $F$, (\ref{E1}) goes back to the stochastic Hamiltonian system
\beq\label{HS} \beg{cases} \d X_t= \nn H(X_t, \cdot)(Y_t)\d t,\\
\d Y_t =-\big\{\nn H(\cdot,Y_t)(X_t)+F(X_t,Y_t)\nn H(X_t,\cdot)(Y_t)\big\}\d t +\d B_t\end{cases}\end{equation} with Hamiltonian function $H$. See e.g. \cite{So} for the physical background and applications in mechanics of the model, and see \cite{T} for exponential convergence of the system to the invariant probability measure. In particular, if $H(x,y)= V(x)+\ff 1 2|y|^2$ and $F\equiv c$ for some constant $c$, (\ref{HS}) is associated to the $``$kinetic Fokker-Planck equation'' in PDE, see e.g. \cite{V} where the hypocoercivity and related regularization estimates w.r.t. the invariant probability measure are studied; and is known as $``$stochastic damping Hamiltonian system'' in probability theory, see e.g. \cite{BCG, Wu} where some long time behaviors of the system have been investigated.

 Following the line of two recent papers \cite {Z,GW} where Bismut formula and Harnack inequalities
are derived for $P_t$ associated to (\ref{E1}) with $Z^{(1)}(x,y)=Ay$ for some $m\times d$-matrix $A$, we aim to derive explicit point-wise derivative estimates of $P_t$ for  more general settings where $Z^{(1)}(x,y)$ might be non-linear and depend on both variables $x$ and $y$, so that some typical examples for the physical model (\ref{HS}) are covered (see Example 4.1 below).

To compare the present equation with those investigated in \cite{Z,GW} where $Z^{(1)}$ is linear, let us recall some simple notations. Firstly, we write the gradient operator on $\R^{m+d}$ as  $\nn =(\nn^{(1)}, \nn^{(2)})$, where $\nn^{(1)}$ and $\nn^{(2)}$ stand for the gradient operators for the first and the second components respectively, so that $\nn f: \R^{m+d}\to\R^{m+d}$ for a differentiable function $f$ on $\R^{m+d}$. Next, for a smooth function $\xi=(\xi_1,\cdots, \xi_k): \R^{m+d}\to \R^k$, let
$$\nn\xi = \left(\beg{matrix} \nn \xi_1\\
\vdots\\
\nn\xi_k\end{matrix}\right),\ \ \nn^{(i)} \xi = \left(\beg{matrix} \nn^{(i)} \xi_1\\
\vdots\\
\nn^{(i)}\xi_k\end{matrix}\right),\ \ \ i=1,2.$$ Then $\nn \xi, \nn^{(1)}\xi, \nn^{(2)}\xi$ are matrix-valued functions of
orders $k\times (m+d), k\times m, k\times d$ respectively. Moreover, for an $l\times k$-matrix $M=(M_{ij})_{1\le i\le l, 1\le j\le k}$
and $v=(v_i)_{1\le i\le k}\in \R^k$, let $Mv\in\R^l$ with $(Mv)_i= \sum_{j=1}^k M_{ij}v_j, \ 1\le i\le l.$  Finally, we will use $\|\cdot\|$  to
denote the operator norm for linear operators, for instance, $\|M\|=\sup_{|v|=1}|Mv|$.

When $Z^{(1)}(x^{(1)},x^{(2)})$ depends only on $x^{(2)}$ and $\nn^{(2)}Z^{(1)}$ is a constant matrix with rank $m$,
then equation (\ref{E1})  reduces back to the one studied in \cite{GW} (and also in \cite{Z} for $m=d$).   In this case we are able to construct very explicit successful couplings with control,
which imply the desired derivative formula and Harnack inequalities as in the elliptic case.
But when $Z^{(1)}$ is non-linear, it seems very hard to construct such couplings.
The idea of this paper is to split $Z^{(1)}$ into a linear term and a non-linear term, and to derive
an explicit derivative formula by controlling  the non-linear part using  the linear part  in a reasonable way. More precisely, let
$$ \nn^{(2) } Z^{(1)}= B_0+B,$$  where $B_0$ is a constant $m\times d$-matrix. We will be able to establish derivative formulae for $P_t$ provided $B$ is dominated by $B_0$ in the  sense  that
\beq\label{B} \<BB_0^* a,a\>\ge -\vv |B_0^*a|^2,\ \ \forall a\in \R^m\end{equation}   holds for some constant $\vv\in [0,1).$

To state our main result, we  first briefly recall the integration by parts formula for the Brownian motion.
Let $T>0$ be fixed. For an Hilbert space $H$, let
$$\H(H)=\bigg\{h\in C([0,T];H):\ h_0=0, \|h\|_{\H(H)}^2:=\int_0^T |\dot h_t|_H^2 \d t<\infty\bigg\}$$ be the Cameron-Martin space over $H$. Let $\H=\H(\R^d)$ and, without confusion in the context,  simply denote  $\|\cdot\|_\H=\|\cdot\|_{\H(H)}$   for any Hilbert space $H$.

Let $\mu$ be the distribution of $\{B_t\}_{t\in [0,T]}$, which is a probability measure (i.e. Wiener measure) on the path space  $\Omega=C([0,T];\R^d)$. The probability space $(\OO,\mu)$ is endowed with the natural filtration of the coordinate process $B_t(w):=w_t, t\in [0,T].$
A function $F\in L^2(\Omega;\mu)$ is called differentiable if for any $h\in \H$, the directional derivative
$$D_h F:= \lim_{\vv\to 0} \ff{F(\cdot+\vv h)-F(\cdot)}{\vv}$$ exists in $L^2(\Omega;\mu)$. If the map $\H\ni h\mapsto D_h F\in L^2(\Omega;\mu)$ is bounded,
then there exists a unique   $DF\in L^2(\Omega\to\H;\mu)$ such that $\<DF, h\>_\H= D_h F$ holds
in $L^2(\Omega;\mu)$ for all $h\in\H$. In this case we write $F\in \D(D)$ and call $DF$
the Malliavin gradient of $F$. It is well known that $(D,\D(D))$ is a closed operator in $L^2(\Omega;\mu)$,
whose adjoint operator $(\dd,\D(\dd))$ is called the divergence
operator. That is,

\beq\label{INT}\E(D_h F)= \int_\Omega D_h F\d\mu= \int_\Omega F
\dd(h)\d\mu= \E(F\dd(h)),\ \ \ F\in \D(D), h\in
\D(\dd).\end{equation}

\

For any $s\ge 0,$ let $\{K(t,s)\}_{t\ge s}$ solve the following random ODE on $\R^m\otimes\R^m$:
\begin{align}
\ff{\d}{\d t}K(t,s)=  (\nn^{(1)}Z^{(1)})(X_t) K(t,s),\ \ \ K(s,s)=I_{m\times m}.\label{Eq1}
\end{align}
We  assume
\paragraph{(H)} The matrix $\sigma\in\R^d\otimes\R^d$ is invertible, and there exists
$W\in C^2(\R^{m+d})$ with $W\ge 1$ and $\lim\limits_{|x|\to\infty} W(x)=\infty$ such that
for some constants $C,l_2\ge 0$ and $l_1\in [0,1],$
\beg{enumerate}
\item[(H1)] $LW\le CW,\ |\nn^{(2)}W|^2\le CW$, where $L=\frac{1}{2}\mathrm{Tr}(\sigma\sigma^*\nabla^{(2)}\nabla^{(2)})+Z\cdot\nabla$;
\item[(H2)] $\|\nn Z\|\le CW^{l_1},\ \ \|\nn^2 Z\|\le CW^{l_2}$.
\end{enumerate}

\

For any $v=(v^{(1)},v^{(2)})\in \R^{m+d}$ with $|v|=1$, we aim to search for $h=h(v)\in\D(\dd)$ such that
 \beq\label{BS} \nn_v P_T f(x)= \E \big[f(X_T(x)) \dd(h)\big],\ \ \ f\in C_b^1(\R^{m+d})\end{equation}
 holds. To construct $h$, for an $\H$-valued random variable $\aa=(\aa_s)_{s\in [0,T]}$, let
 \beg{equation}\label{B0}\beg{split}& g_t= K(t,0)v^{(1)} +\int_0^tK(t,s)\nn^{(2)}Z^{(1)}(X_s(x))\aa_s\d s,\\
& h_t=\int_0^t\si^{-1}\big(\nn Z^{(2)}(X_s(x)) (g_s,\aa_s) -\dot\aa_s\big)\d s,\ \ t\in [0,T].\end{split}\end{equation}
We will show that $h$ satisfies (\ref{BS}) provided it is in $\D(\dd)$ and $\aa_0=v^{(2)},\aa_T=0,g_T=0$, see
Theorem \ref{T1.1} below for details. In particular, it is the case for $\aa_s$ given in the following result.

\beg{thm} \label{T1.2} Assume {\bf (H)} and let $\nn^{(2)}Z^{(1)}=B_0+B$ for some constant matrix $B_0$ such that
$(\ref{B})$  holds for some constant $\vv\in [0,1).$
If there exist an increasing fcuntion $\xi\in C([0,T])$ and $\phi\in C^1([0,T])$ with $\xi(t)>0$ for $t\in (0,T]$, $\phi(0)=\phi(T)=0$ and $\phi(t)>0$ for $t\in (0,T)$
such that
\beq\label{B2} \int_0^t \phi(s) K(T,s)B_0B_0^*K(T,s)^*\d s\ge \xi(t) I_{m\times m},\ \ t\in (0,T].\end{equation}
Then \beg{enumerate}\item[$(1)$] $Q_t:=\int_0^t \phi(s)K(T,s)\nn^{(2)}Z^{(1)}(X_s) B_0^* K(T,s)^* \d s $ is invertible for $t\in (0,T]$ with
\beq\label{Q}\|Q_t^{-1}\|\le\ff 1 { (1-\vv)\xi(t)},\ \ t\in [0,T].\end{equation}
\item[$(2)$] Let    $h$ be determined by $(\ref{B0})$ for
\beg{equation}\label{aa}\beg{split} \aa_t := &\ff{T-t} T v^{(2)} -\phi(t) B_0^* K(T,t)^*Q_T^{-1}\int_0^T\ff{T-s}TK(T,s) \nn^{(2)}Z^{(1)}(X_s) v^{(2)} \d s\\
&-\ff{\phi(t)B_0^*K(T,t)^*}{\int_0^T \xi(s)^2 \d s} \int_t^T  \xi(s)^2 Q_s^{-1} K(T,0)v^{(1)}\d s.
\end{split}\end{equation}
Then for any $p\ge 2$, there exists a constant $T_p\in (0,\infty)$ if $l_1=1$ and $T_p=\infty$ if $l_1<1$,
such that for any $T\in (0,T_p)$, $(\ref{BS})$ holds  with $\E|\dd(h)|^p<\infty$.
\item[$(3)$] For any $p>1$ there exist    constants $c_1(p),c_2(p)\ge 0$, where $c_2(p)=0$ if $l_1=l_2=0$,  such that
\beq\label{GG}|\nn P_T f|\le c_1(p)(P_T| f|^p)^{1/p}  \ff{\ss{T\land 1}\{(T\land 1)^2 +\xi(T\land 1)\}\e^{c_2(p)W}}{\int_0^{T\land 1} \xi(s)^2\d s}\end{equation} holds for all $ T>0$ and $f\in \B_b(\R^{m+d}).$\end{enumerate}  \end{thm}

The remainder of the paper is organized as follows.  In Section 2 we present a general result on
the derivative formula by using Malliavin calculus, from which we are able to prove Theorem \ref{T1.2}  in Section 3.  In Section  4  we will verify    (\ref{B2}) for the following two cases respectively:

\beg{enumerate}\item[(I)] $\nn^{(1)}Z^{(1)}$ is non-constant but $\text{Rank}[B_0]=m.$
\item[(II)] $A:= \nn^{(1)}Z^{(1)}$ is constant such that $\text{Rank}[B_0, AB_0,\cdots, A^kB_0]=m$ holds for some $0\le k\le m-1.$
\end{enumerate}
In both cases the $L^p$-gradient estimate (\ref{GG}) is derived with specific $\xi$, while in Case (II)
the Harnack inequality introduced in \cite{W97} is established provided
$\nn Z^{(1)}$ is constant, which extends the corresponding Harnack inequality obtained in
\cite{GW} for $\nn^{(1)}Z^{(1)}=0$ and $\nn^{(2)}Z^{(1)}$ is constant with rank $m$. This type of Harnack inequality has been applied in  the study of heat kernel estimates and contractivity properties of Markov semigroups, see e.g. \cite{GW} and references therein.

\section{A General Result}

In this section we will make use of the following assumption.
\paragraph{{\bf (H$'$)}}   The function
$$U(x):= \E\exp\bigg[2\int_0^T\|\nn Z(X_t(x))\|\d t\bigg] ,\ \ x\in \R^{m+d}$$ is locally bounded.

\beg{thm}\label{T1.1} Assume  {\bf{\bf (H$'$)}} for some $T>0$. For   $v=(v^{(1)}, v^{(2)})\in \R^{m+d}$,  let  $(\aa_s)_{0\le s\le T}$ be   an $\H$-valued random variable such that
 $\aa_0=v^{(2)}$ and $\aa_T=0,$ and let $g_t$ and $h_t$ be given in $(\ref{B0})$.
 If $g_T=0$ and $h\in \D(\dd)$, then $(\ref{BS})$ holds.
\end{thm}

 \beg{proof} For simplicity, we will drop the initial data of the solution by writing $X_t(x)=X_t$.
 By {\bf {\bf (H$'$)}} and (\ref{E2}) we have
 $X_t\in \D(D)$,  and due to the chain rule and the definition of $h_t$,
  \beg{equation}\label{P1}\beg{split} D_hX_t&= \int_0^t \nn Z(X_s) D_h X_s\d s +\int_0^t (0, \si \dot h_s)\d s\\
  &=  (0, v^{(2)}-\aa_t)+ \int_0^t \nn Z(X_s) D_h X_s\d s +\int_0^t  \left(0,  \nn Z^{(2)}(X_s)(g_s,\aa_s)\right)\d s\end{split}\end{equation}
holds for $t\in [0,T].$  Next, it is easy to see that
 $$g_t= v^{(1)} + \int_0^t  \nn Z^{(1)}(X_s) (g_s,\aa_s) \d s,\ \ \ t\in [0,T].$$ Combining this with (\ref{P1}) we
 obtain
 $$D_h X_t +(g_t,\aa_t) = v +\int_0^t \nn Z(X_s) \{D_h X_s +(g_s,\aa_s)\} \d s,\ \ t\in [0,T].$$ On the other hand, the directional derivative process
 $$\nn_v X_t:=\lim_{\vv\to 0} \ff{X_t(x+\vv v)-X_t(x)}{\vv}$$ satisfies the same equation, i.e.
\beq\label{**0}\nn_v X_t = v+ \int_0^t \nn Z(X_s) \nn_vX_s\d s,\ \ \ t\in [0,T].\end{equation} Thus, by the uniqueness of the ODE we conclude that
 $$D_h X_t+(g_t,\aa_t)= \nn_v X_t,\ \ \ t\in [0,T].$$ In particular, since $(g_T,\aa_T)=0$, we have
 \begin{align}
 D_h X_T= \nn_v X_T\label{EQ}
 \end{align}
 and due to {\bf {\bf (H$'$)}} and (\ref{**0}),
 \beq\label{**1} \E|D_h X_T|^2=\E|\nn_vX_T|^2 \le |v|^2\E\exp\bigg[2\int_0^T\|\nn Z\|(X_s)\d s\bigg].\end{equation} Combining this with (\ref{INT}) and letting $f\in C_b^1(\R^{m+d})$, we are able to adopt the dominated convergence theorem to obtain
$$\nn_v P_T f = \E   \<\nn f(X_T),\nn_v X_T\>= \E\<\nn f(X_T), D_h X_T\> = \E D_h f(X_T)= \E[f(X_T)\dd(h)].$$
\end{proof}

\beg{rem}
Using the same argument as above, we also have the following derivative formula:
\begin{align}
\E\nabla_v f(X_T)=\E\left(f(X_T)\sum_{i,k}\Big[\delta(h(e_k))(\nabla X_T)^{-1}_{ki}
-D_{h(e_k)}(\nabla X_T)^{-1}_{ki}\Big]v^i\right),\label{Ep4}
\end{align}
where $(e_j)$ is the canonical basis of $\R^{m+d}$, and $h(e_j)$ is defined by (\ref{B0}) with $v=e_j$.
In fact, since
$$
\sum_{k}(\partial_k X^j_T)(\nabla X_T)^{-1}_{ki}=1_{i=j}
$$
and by (\ref{EQ})
$$
D_{h(e_k)}X^j_T=\nabla_{e_k}X^j_T=\partial_k X^j_T,
$$
we have
\begin{align*}
\nabla_v f(X_T)=\sum_i(\partial_i f)(X_T)v^i
&=\sum_{i,j,k}(\partial_j f)(X_T)(\partial_k X^j_T)(\nabla X_T)^{-1}_{ki}v^i\\
&=\sum_{i,j,k}(\partial_j f)(X_T)(D_{h(e_k)}X^j_T)(\nabla X_T)^{-1}_{ki}v^i\\
&=\sum_{i,k}\{D_{h(e_k)} f(X_T)\}(\nabla X_T)^{-1}_{ki}v^i,
\end{align*}
which implies (\ref{Ep4}) by the integration by parts formula.
\end{rem}
\beg{rem}
For the higher order derivative formula, under further regularity assumptions,
for any $v_1,\cdots,v_j\in\R^{m+d}$ and $f\in C^1_b(\R^{m+d})$, we have
\begin{align}
\<\nabla^j\E f(X_T(x)),v_1\otimes\cdots\otimes v_j\>=
\E\left[f(X_T(x)) J_j(T,v_1,\cdots,v_j)\right],\label{FOR1}
\end{align}
where $J_1(v):=\delta(h(v))$ and
\begin{align*}
J_j(v_1,\cdots,v_j)&:=J_{j-1}(v_1,\cdots,v_{j-1})\delta(h(v_j))
+\nabla_{v_j} J_{j-1}(v_1,\cdots,v_{j-1})\\
&\quad-D_{h(v_j)}J_{j-1}(v_1,\cdots,v_{j-1}),
\end{align*}
where $h(v)$ is defined by (\ref{B0}). In fact, as in the proof of Theorem \ref{T1.1}, we have
\begin{align*}
&\<\nabla^2\E f(X_{T}),v_1\otimes v_2\>=\nabla_{v_2}\nabla_{v_1}\E f(X_{T})
=\nabla_{v_2}\E[f(X_{T})\delta(h(v_1))]\\
&\quad=\E\left[(\nabla f)(X_{T})\cdot\nabla_{v_2} X_{T}\cdot\delta(h(v_1))\right]
+\E\left[f(X_{T})\nabla_{v_2}\delta(h(v_1))\right]\\
&\quad=\E\left[(\nabla f)(X_{T})\cdot D_{h(v_2)}X_{T}\cdot\delta(h^{v_1})\right]
+\E\left[f(X_{T}(x))\nabla_{v_2}\delta(h(v_1))\right]\\
&\quad=\E\left[D_{h(v_2)}[f(X_{T})]\delta(h(v_1))\right]
+\E\left[f(X_{T}(x))\nabla_{v_2}\delta(h(v_1))\right]\\
&\quad=\E\left[f(X_{T}(x))\big[\delta(h(v_1))\delta(h^{v_2}_T)
-D_{h(v_2)}\delta(h(v_1))+\nabla_{v_2}\delta(h(v_1))\big]\right].
\end{align*}
The higher derivatives can be obtained by induction.
\end{rem}

\section{Proof of Theorem \ref{T1.2}}
The idea of the proof is to apply Theorem \ref{T1.1} for the given process $\aa_s$.   Obviously, (H1) implies that for any $l\geq 1$, there exists a constant $C_l$ such that
$LW^l\le C_l W^l$, so that $\E W(X_t(x))^l\le \e^{C_lt}W(x)^l$ and thus, the process is non-explosive; while (H2)
imply that $\|\nn Z\| +\|\nn^2 Z \|\le CW^{l_1\lor l_2}$ holds for some $C>0$, so that
\beq\label{G0} \E\left( \big(\|\nn Z\|^p+\|\nn^2 Z\|^p\big)(X_t) \right)\le \e^{c(p)t} W^{p(l_1\lor l_2)},\ \ t\ge 0\end{equation} holds for
any $p\geq 1$ with some constant $c(p)>0.$ The following lemma ensures that {\bf(H)} implies {\bf {\bf (H$'$)}}
for all $T>0$ if $l_1<1$ and for small $T>0$ if $l_1=1$.

\beg{lem}\label{L1} If {\rm (H1)} holds, then for any $T>0$,
$$\E\exp\bigg[\ff 2 {T^2C\|\si\|^2\e^{4+2CT}}\int_0^T W(X_t)\d t\bigg]\le \exp\bigg[\ff{2W}{TC\|\si\|^2\e^{2+CT}}\bigg].$$
Consequently, {\rm (H2)} imply that $U:= \E\exp[2\int_0^T\|\nn Z\|(X_t)\d t]$ is locally bounded on $\R^{m+d}$ if either $l_1<1$ or
$l_1=1$ but $T^2C^2\|\si\|^2\e^{4+2CT}\le 1.$
\end{lem}

\beg{proof} It suffices to prove the first assertion. By the It\^o formula and (H1), we have
$$\d W(X_t)= \<\nn^{(2)} W(X_t), \si\d B_t\> + LW(X_t)\d t\le  \<\nn^{(2)} W(X_t), \si\d B_t\> + C W(X_t)\d t.$$ So, for $t\in [0,T]$,
$$ \d\big\{\e^{-(C+2/T)t}W(X_t)\big\}\le \e^{-(C+2/T)t}\<\nn^{(2)} W(X_t), \si\d B_t\>- \ff 2 T \e^{-CT-2} W(X_t)\d t.$$ Thus, letting
$\tau_n=\inf\{t\ge 0:W(X_t)\ge n\}$, for any $n\ge 1$ and  $\ll>0$ we have
\beg{equation*}\beg{split} &\E\exp\bigg[\ff{2\ll}{T\e^{CT+2}}\int_0^{T\land\tau_n} W(X_t)\d t\bigg]\\
&\le \e^{\ll W}\E \exp\bigg[\ll\int_0^{T\land\tau_n}\e^{-(C+2/T)t}\<\nn^{(2)} W(X_t), \si\d B_t\>\bigg]\\
&\le \e^{\ll W} \bigg(\E\exp\bigg[2\ll^2C\|\si\|^2 \int_0^{T \land\tau_n}W(X_t)\d t\bigg]\bigg)^{1/2},\end{split}\end{equation*}
where the second inequality is due to the exponential martingale and (H1). By taking
$$\ll= \ff 1 {TC\|\si\|^2 \e^{CT+2}},$$ we arrive at
$$\E\exp\bigg[\ff 2 {T^2C\|\si\|^2\e^{4+2CT}}\int_0^{T\land\tau_n} W(X_t)\d t\bigg]\le \exp\bigg[\ff{2W}{TC\|\si\|^2\e^{2+CT}}\bigg].$$
This completes the proof by letting $n\to\infty.$
\end{proof}
To ensure that $\E|\dd(h)|^p<\infty$, we need the following two lemmas.
\beg{lem}\label{L2} Assume {\bf (H)}. Then  there exists a constant $c>0$ such that
\beq\label{G1}\|DX_t\|_{\H}\le \sqrt{t}\|\si\|\e^{c\int_0^tW^{l_1}(X_s)\d s}, t\geq 0.\end{equation}
Consequently, if $l_1<1$, then for any $p\geq 1$,
$$\E\left(\sup_{t\in [0,T]} \|DX_t\|_{\H}^p\right)<\infty,\ \ T\geq 0;$$ and if $l_1=1$,
then for any $p\geq1$ there exists a constant $T_p>0$ such that
$$\E \left(\sup_{t\in [0,T]}\|DX_t\|_{\H}^p\right)<\infty,\ \ T\in (0,T_p).$$ \end{lem}

\beg{proof} Due to Lemma \ref{L1}, it suffices to prove (\ref{G1}). From (\ref{E2}) we see that for any $h\in\H$, $D_hX_t$ solves the following   random ODE:
$$ D_hX_t = \int^t_0(\nn Z)(X_s) D_hX_s\d s +(0,\si h(t)).$$ Combining this with (H2) and   $|h(t)|\le \ss{t}\,\|h\|_\H$, we obtain
$$|D_hX_t| \le C\int^t_0W^{l_1}(X_s) |D_hX_s|\d s +\sqrt{t}\|\si\|\cdot \|h\|_\H,\ \ h\in\H.$$Therefore,
$$\|DX_t\|_{\H} \le C\int^t_0W^{l_1}(X_s) \|DX_s\|_{\H}\d s +\sqrt{t}\|\si\|.$$
This implies (\ref{G1}) by Gronwall's inequality.\end{proof}

\beg{lem}\label{L3} Assume {\bf(H)}. Then for any $s\in[0,T]$,
\beq\label{Es1}
\|K(T,s)\|\leq C\e^{C\int^T_s W^{l_1}(X_r)\d r},\ \|\partial_sK(T,s)\|\leq CW^{l_1}(X_s)\e^{C\int^T_s W^{l_1}(X_r)\d r},
\end{equation}
and
\beq\label{G2}
\|DK(T,s)\|_{\H}\le C\e^{C\int^T_s W^{l_1}(X_r)\d r} \int_s^T W^{l_2}(X_r)\|DX_r\|_{\H}\d r.
\end{equation}
Consequently, for any $p>1$ there exists $T_p\in (0,\infty)$ if $l_1=1$ and $T_p=\infty$ if $l_1<1$ such that
$$\E\left(\sup_{t\in [0,T]}\|DK(T,t)\|_{\H}^p\right)<\infty,\ \ T\in (0,T_p).$$\end{lem}
\beg{proof}
By Lemma \ref{L2} and $\sup_{t\in [0,T]} \E W^l (X_t)<\infty$ for any $l>0$ as observed in the beginning of this section, it suffices to
prove (\ref{Es1}) and (\ref{G2}).
First of all, by (\ref{Eq1}) and (H2), we  have
\begin{align*}
\|K(t,s)\|\leq 1+\int^t_s\|\nn^{(1)}Z^{(1)}(X_r)\|~\|K(r,s)\|\d r\leq 1+C\int^t_sW^{l_1}(X_r)\|K(r,s)\|\d r.
\end{align*}
which yields the first estimate in  (\ref{Es1})  by Gronwall's inequality.
Moreover, noticing that
$$
\partial_s K(t,s)=\int^t_s(\nn^{(1)}Z^{(1)})(X_r) \partial_s K(r,s)\d r-(\nn^{(1)}Z^{(1)})(X_s),
$$
by (H2) we have
$$
\|\partial_s K(t,s)\|\leq C\int^t_sW^{l_1}(X_r)\|\partial_s K(r,s)\|\d r+CW^{l_1}(X_s).
$$
The second estimate in  (\ref{Es1}) follows.  As for (\ref{G2}), since
$$\ff{\d}{\d t} DK(t,s)= (\nn_{DX_t}\nn^{(1)}Z^{(1)})(X_t) K(t,s) +(\nn^{(1)}Z^{(1)})(X_t) D K(t,s),$$
with $DK(s,s)=0$, it follows from (H2) and (\ref{Es1}) that
\beg{equation*}\beg{split} \|DK(t,s)\|_{\H}&\le
\int^t_s\|\nn\nn^{(1)}Z^{(1)}(X_r)\|~\|DX_r\|_{\H}\|K(r,s)\|\d r\\
&\quad+\int^t_s\|\nn^{(1)}Z^{(1)}(X_r)\|~\|DK(r,s)\|_{\H}\d r \\
&  \le C\e^{C\int^T_s W^{l_1}(X_r)\d r} \int^t_s W^{l_2}(X_r)\|DX_r\|_{\H} \d r\\
&\quad+C\int^t_sW^{l_1}(X_r)\|DK(r,s)\|_{\H}\d r.\end{split}\end{equation*}
This implies (\ref{G2}).\end{proof}

\beg{proof}[Proof of Theorem \ref{T1.2}]
(1) Let $a\in \R^m$. By (\ref{B}), (\ref{B2}) and $\nn^{(2)}Z^{(1)}=B_0+B$  we have

\beg{equation*}\beg{split} \<Q_t a,a\> &= \int_0^t \phi(s) \Big(\<K(T,s)B_0B_0^*K(T,s)^*a,a\>+\<K(T,s)B(X_s)B_0^*K(T,s)^*a,a\>\Big)\d s\\
&\ge (1-\vv)\int_0^t \phi(s) |B_0^*K(T,s)^*a|^2\d s \ge (1-\vv)\xi(t)|a|^2.\end{split}\end{equation*} This implies that $Q_t$ is invertible and
(\ref{Q}) holds.

(2)  According to Lemma \ref{L1}, {\bf (H)} implies {\bf  (H$'$)} for all $T>0$ if $l_1<1$ and
for small $T>0$ if $l_1=1.$ Next, we intend to prove that $h\in \D(\dd)$ and $\E |\dd(h)|^p<\infty$ for small
$T>0$ if $l_1=1$ and for all $T>0$ if $l_1<1.$
Indeed, by Lemmas \ref{L2}, \ref{L3}, (\ref{G0}), and the fact that
$$
DQ_t^{-1}=-Q_t^{-1}(DQ_t) Q_t^{-1},
$$
there exists $T_p>0$ if $l_1=1$ and $T_p=\infty$ if $l_1<1$ such that
$$
\sup_{t\in[0,T]}\E |DQ_t|^p<+\infty,\ \ T\in (0,T_p),
$$
and by (\ref{Q}),
\beq\label{AG01}
\Big(\E\|D Q_t^{-1}\|_{\H}^p\Big)^{1/p}\leq \frac{\big(\E|DQ_t|^p\big)^{1/p}}{[(1-\epsilon)\xi(t)]^2},\ \ t\in(0,T],
\end{equation}
\beq\label{AG0}
\sup_{t\in [0,T]}\Big(\E\|D\aa_t\|_{\H}^p+\E\|Dg_t\|_{\H}^p\Big)^{1/p}<\infty,\ \ T\in (0,T_p).
\end{equation}
Since
\beq\label{AG}\beg{split} & \dot h_t = \si^{-1}\big\{(\nn Z^{(2)})(X_t)(g_t,\aa_t)-\dot\aa_t\big\},\\
&\|D\dot h_t\|_{\H}\le \|\si^{-1}\|\big\{\|\nn^2Z^{(2)}(X_t)\|~\|DX_t\|_{\H}~ |(g_t,\aa_t)|\\
&\qquad\qquad\quad+\|\nn Z^{(2)}(X_t)\|~ \|(Dg_t,D \aa_t)\|_{\H}+ \|D\dot\aa_t\|_{\H}\big\}.\end{split}\end{equation}
we conclude  from (H2),  (\ref{G0}) and  (\ref{AG0}) that
$$\E\bigg(\int_0^T \|D\dot h_t\|_\H^2\d t \bigg)^{p/2}+ \E\|h\|_\H^p<\infty,\ \ T\in (0,T_p).$$   Therefore,
according to e.g.  \cite[Proposition 1.5.8]{N},  we have $h\in\D(\dd)$ and $\E|\dd(h)|^p<\infty$ provided $T\in (0,T_p).$

Now, to prove (\ref{BS}),  it remains to verify the required conditions of Theorem \ref{T1.1} for $\aa_t$   given by (\ref{aa}). Since $\phi(0)=\phi(T)=0$, we have
 $\aa_0=v^{(2)}$ and $\aa_T=0.$ Moreover,  noting that
 \beg{equation*}\beg{split}I_1&:= \ff 1 {\int_0^T \xi(t)^2\d t}
 \int_0^T \phi(t)K(T,t)\nn^{(2)}Z^{(1)}(X_t) B_0^*K(T,t)^*\d t\int_t^T\xi(s)^2Q_s^{-1}K(T,0)v^{(1)}\d s\\
 &=\ff 1 {\int_0^T \xi(t)^2\d t}\int_0^T \dot Q_t\d t \int_t^T\xi(s)^2Q_s^{-1}K(T,0)v^{(1)}\d s\\
 &= \ff 1 {\int_0^T \xi(t)^2\d t}\int_0^T  \xi(t)^2Q_tQ_t^{-1}K(T,0)v^{(1)}\d t
 =K(T,0)v^{(1)}\end{split}\end{equation*} and

 \beg{equation*}\beg{split}  I_2&:=\bigg(\int_0^T\phi(t)K(T,t)\nn^{(2)}Z^{(1)}(X_t)B_0^*K(T,t)^*\d t\bigg)Q_T^{-1}\int_0^T\ff{T-s}TK(T,s)\nn^{(2)}Z^{(1)}(X_s)v^{(2)}\d s\\
 &=Q_TQ_T^{-1} \int_0^T\ff{T-s}TK(T,s)\nn^{(2)}Z^{(1)}(X_s)v^{(2)}\d s= \int_0^T\ff{T-s}TK(T,s)\nn^{(2)}Z^{(1)}(X_s)v^{(2)}\d s,
 \end{split}\end{equation*}
we obtain by (\ref{aa})
\beg{equation*}\beg{split} g_T&= K(T,0)v^{(1)}+\int_0^TK(T,t)\nn^{(2)}Z^{(1)} (X_t)\aa_t\d t \\
&=K(T,0)v^{(1)} -I_1 +  \int_0^T\ff {T-t}T K(T,t) \nn^{(2)}Z^{(1)}(X_t)v^{(2)}\d t- I_2=0.\end{split}\end{equation*}

(3) By an approximation argument, it suffices to prove the desired gradient estimate for $f\in C_b^1(\R^{m+d}).$
Moreover, by the semigroup property and the Jensen inequality, we only have to prove for  $p\in(1,2]$ and $T\in (0, T_p\land 1).$ In this case we obtain from
(\ref{BS}) that $$|\nn P_T f|\le (P_T |f|^p)^{1/p} (\E|\dd(h)|^q)^{1/q},$$ where
$q:=\ff p{p-1}\ge 2.$  Therefore, it remains to find  constants $c_1,c_2\ge 0$, where $c_2=0$ if $l_1=l_2=0$,
such that
\beq\label{WW} (\E |\dd(h)|^q)^{1/q} \le \ff{c_1\ss T(T^2+\xi(T))\e^{c_2W}}{ \int_0^T\xi(s)^2\d s }.\end{equation} To this end, we take $\phi(t)= \ff{t(T-t)}{T^2}$ such that $0\le \phi\le 1$
and $|\dot \phi(t)|\le \ff 1 T$ for $t\in [0,T].$ Since $\xi$ is increasing, by (\ref{Es1}) and (\ref{B2}),
we have for some constant $C>0$,
$$
\int^t_0\xi(s)^2\d s\leq \xi(t)^2\leq Ct^2,\ \ t\in[0,1].
$$
 Thus, by Lemmas \ref{L1}, \ref{L2}, \ref{L3} and  (\ref{G0}), it is easy to see that for any $\theta\ge 2$ there exist  constants  $c_1,c_2\ge 0$, where
 $c_2=0$ if $l_1=l_2=0$,  such that for all $0<t\leq T\leq T_p\wedge 1$,
\beg{equation*}\beg{split}& \big(\E\|DX_t\|_{\H}^\theta\Big)^{1/\theta}
 \le c_1\sqrt{T} \e^{c_2W},\ \ \big(\E\|DK(T,t)\|_{\H}^\theta\big)^{1/\theta}\le c_1 T^{3/2} \e^{c_2W}\\
   &( \E \|DQ_t^{-1}\|_{\H}^\theta)^{1/\theta}\le \big\{\E(\|Q_t^{-1}\|\|DQ_t\|_{\H}\|Q_t^{-1}\|)^\theta\big\}^{1/\theta}
 \le \ff{c_1t\sqrt{T}}{\xi(t)^2}\e^{c_2 W},\\
 &(\E\|D\aa_t\|_{\H}^\theta)^{1/\theta} \le\ff{c_1T^{5/2}\e^{c_2W}}{\int_0^T \xi(s)^2\d s},\ \
 (\E\|Dg_t\|_{\H}^\theta  )^{1/\theta}\le \ff{c_1T^{7/2}\e^{c_2W}}{\int_0^T \xi(s)^2\d s}, \\
& (\E\|D\dot\aa_t\|_{\H}^\theta)^{1/\theta}  \le \ff{c_1 T^{3/2}\e^{c_2W}}{\int_0^T \xi(s)^2\d s},\ \
(\E |\dot h_t|^\theta)^{1/\theta} \le \ff{c_1 \xi(T)\e^{c_2W}}{\int_0^T \xi(s)^2\d s}.\end{split}\end{equation*}
Combining these with (\ref{AG}), (H2)  and (\ref{G0}),  we obtain

\beg{equation*}\beg{split}\|h\|_{\mathbb D^{1,q}}&:=\big(\E\|Dh\|_{\H\otimes\H}^q\big)^{1/q} + \|\E h\|_\H\\
&\le \ss T \bigg\{\E\bigg(\ff 1 T\int_0^T \|D\dot h_t\|_{\H}^2\d t\bigg)^{q/2}\bigg\}^{1/q} +\E\|h\|_\H\\
&\le \ss T \bigg(\ff 1 T \int_0^T \E \|D \dot h_t\|_{\H}^q\d t\bigg)^{1/q} +\bigg(\E\int_0^T|\dot h_t|^2\d t\bigg)^{1/2}\\
&\le \ff{c_1\ss T(T^{3/2}+\xi(T))\e^{c_2W}}{\int_0^T\xi(s)^2\d s}.\end{split}\end{equation*}
This implies (\ref{WW})  since $\dd: \mathbb D^{1,q}\to L^q$ is bounded, see e.g.  Proposition 1.5.8 in \cite{N}.
 \end{proof}
\section{Two Specific Cases}

As indicated in the end of Section 1, we intend to apply Theorem \ref{T1.2} to Case (I) and Case (II) respectively with concrete functions $\xi$ satisfying
(\ref{B2}).
\subsection{Case (I): Rank$[B_0]=m$}

\beg{thm} \label{T4.1} Assume {\bf (H)} and $(\ref{B})$ for some $\vv\in [0,1).$ If {\rm Rank}$[B_0]=m$, then there exist constants $c_1,c_2>0$ such that $(\ref{B2})$ holds for
$$\xi(t)= c_1 \int_0^t \phi(s)\e^{-c_2(T-s)}\d s,\ \  t\in[0,T].$$ Consequently, for any $p>1$ there exist two constants $c_1(p),c_2(p)\ge 0$, where $c_2(p)=0$ if
$l_1=l_2=0$, such that
$$|\nn P_T f|\le \ff{c_1(p) (P_T|f|^p)^{1/p}}{(T\land 1)^{3/2}} \e^{c_2(p) W},\ \ \ T>0.$$\end{thm}

\beg{proof} It is easy to see that the desired gradient estimate follows from (\ref{GG}) for the claimed $\xi$ with $\phi(t)=\ff{t(T-t)}{T^2},$ we only prove the first assertion. Since $\nn^{(1)}Z^{(1)}$ is bounded, there exists a constant $C>0$ such that
$$|K(T,s)^*a|\ge \e^{-C(T-s)}|a|,\ \ \ a\in \R^m.$$ If Rank$[B_0]=m$, then $|B_0^*a|\ge c'|a|$ holds for some constant $c'>0$ and all $a\in \R^m$. Therefore,
$$ M_t:= \int_0^t\phi(s) K(T,s) B_0B_0^* K(T,s)^*\d s $$ satisfies
$$\<M_t a,a\> = \int_0^t \phi(s)|B_0^*K(T,s)^*a|^2\d s\ge c'^2\int_0^t \phi(s) \e^{-2C(T-s)}|a|^2\d s.$$ This completes the proof. \end{proof}

\paragraph{Example 4.1.} Consider the stochastic Hamilton system (\ref{HS}),
where $m=d$ and $\nn^{(2)}Z^{(1)}= \Hess_{H(x,\cdot)}(y)$ is symmetric. If
for some $C>0$,
\beq\label{QQL}
CI_{d\times d}\le \nn^{(2)}Z^{(1)},\  \text{or}\ \nn^{(2)}Z^{(1)}\le -C I_{d\times d}.
\end{equation}Then  we take $B_0=CI_{d\times d}$
if $\nn^{(2)}Z^{(1)}\ge CI_{d\times d}$, while $B_0=-CI_{d\times d}$ if $\nn^{(2)}Z^{(1)}\le -CI_{d\times d}$.
It is trivial to see that Rank$[B_0]=d=m$ and (\ref{B}) holds for $\vv=0.$

A typical choice of $H$ in the physical model such that (\ref{QQL})
holds is that (cf. \cite[Chapter XIII]{So})
$$
H(x,y)=V(x)+\frac{1}{2}\<M(x)y,y\>,
$$
where $M(x)$, called mass matrix of the system, is a $d\times d$-real symmetric, smooth and   positive definite matrix; and
$V(x)$, called potential energy, is a smooth function.
Assume that
$$
M(x)\geq CI_{d\times d},$$ then (\ref{QQL}), and hence (\ref{B}) with $\vv=0$,   holds.
If moreover $F\in C_b^2,  M\in C_b^3,$  and $V\ge 0$ (equivalently, bounded from below since one may add a constant to $H$) such that
 $$
\|\nabla^k V(x)\|\leq C(V(x) +1),\ \ k=2,3,
$$
then Assumption {\bf (H)} holds with
$
W(x,y)=H(x,y) +1
$
and $l_1=l_2=1$.  Therefore, Theorem \ref{T4.1} applies.

\subsection{Case (II): $A:=\nn^{(1)}Z^{(1)}$ is constant}

Throughout this subsection we assume that
\paragraph{(A)} (Kalman condition)  $A:=\nn^{(1)}Z^{(1)}$ is constant and there exists an integer number $0\le k\le m-1$ such that
\beq\label{R} \text{Rank}[B_0, AB_0, \cdots, A^{k}B_0]=m.\end{equation}

\  \newline
When $k=0$, (\ref{R})  means  Rank$[B_0]=m$ which has been
considered in Theorem \ref{T4.1}.

\beg{thm} \label{T4.2} Assume    {\bf (H), (A)} and $(\ref{B})$ for some $\vv\in (0,1).$ Let $\phi(t)=\ff{t(T-t)}{T^2}.$ Then:
\beg{enumerate}\item[$(1)$]  There exist constants $c_1,c_2>0$ such that $(\ref{B2})$ holds for
$$\xi(t)=\ff{ c_1 (t\land 1)^{2(k+1)} }{T\e^{c_2 T}},\ \ \ t\in [0,T].$$
\item[$(2)$] For any $p>1$, there exist two constants $c_1(p),c_2(p)\ge 0$, where $c_2(p)=0$ if
$l_1=l_2=0$, such that
$$|\nn P_T f|\le \ff{c_1(p) (P_T|f|^p)^{1/p}}{(T\land 1)^{(4k-1)\vee 0+3/2}} \e^{c_2(p) W},\ \ \ T>0.$$
\item[$(3)$] If $\nn^{(2)}Z^{(1)}=B_0$ is constant and $l_1<\ff 1 2,$ then there exists a constant $c>0$ such that
\beg{equation*}\beg{split}& |\nn P_Tf|\le \ll\big\{P_Tf\log f-(P_Tf)\log P_T f\big\} \\
&+\ff c  \ll \bigg\{\ff{ l_1W}{(1+\ll^{-1})^2}+\ff{(1+\ll^{-1})^{4l_1/(1-2l_1)}}{(T\land 1)^{(4k+2-2l_1)/(1-2l_1)}}+\ff 1 {(1\land T)^{4k+3}}\bigg\}P_Tf,\  \ll>0, T>0 \end{split}\end{equation*}holds for all
 $f\in \B_b^+(\R^{m+d}),$ the set of   positive functions in $\B_b(\R^{m+d}).$
\item[$(4)$] If $\nn^{(2)}Z^{(1)}=B_0$ is constant and $l_1=\ff 1 2,$ then there exist  constants $c,c'>0$ such that for any $T>0,
 \ll\ge \ff c{(T\land 1)^{2k}} $ and  $f\in \B_b^+(\R^{m+d}),$
$$|\nn P_Tf|\le \ll\big\{P_Tf\log f-(P_Tf)\log P_T f\big\}+ \ff{c'((1\land T)^2W+1)}{\ll (T\land 1)^{4k+3}}P_T f.$$
\end{enumerate}
\end{thm}
\beg{proof} Since (2) is a direct consequence of (\ref{GG}) and (1), we only prove (1), (3) and (4).

(1) Let
$$M_t=\int_0^t \ff{s(T-s)}{T^2} \e^{(T-s)A}B_0 B_0^* \e^{(T-s)A^*} \d s,\ \
U_t= \int_0^t \e^{sA}B_0B_0^*\e^{sA^*}\d s,\ \ t\in [0,T].$$ According to \cite[\S3]{S},
the limit
$$Q:=\lim_{t\to 0} t^{-(2k+1)}\GG_tU_t\GG_t$$ exists and is an invertible matrix, where $(\GG_t)_{t>0}$ is a family of projection matrices. Thus,
  $U_t\ge c (t\land 1)^{2k+1} I_{m\times m}$ holds for some constant $c>0$ and all $t>0.$ Then there exist constants $c_1,c_2>0$ such that for any  $t\in (0,\ff T 2]$,
  $$M_t\ge \ff{t}{4T} \int_{t/2}^t \e^{(T-s)A}B_0B_0^* \e^{(T-s)A^*}\d s
  \ge \ff{t\e^{-2\|A\|T}}{4T} \int_0^{t/2} \e^{sA}B_0B_0^*\e^{sA^*}\d s\ge \ff{c_1 t^{2(k+1)}}{4T\e^{c_2T}}I_{m\times m}$$ holds.
  This proves the first assertion.

(3) By the semigroup property and the Jensen inequality, we assume that $T\in (0,1].$ Let $\nn^{(2)}Z^{(1)}=B_0$ be constant. Then $h$ given in Theorem \ref{T1.2} is adapted such that
$$\dd(h)=\int_0^T \<\dot h_t,\d B_t\> .$$ Moreover, it is easy to see that for $\xi(t)$ given in (1) and $T\in (0,1],$
$$|\dot h_t|\le \ff{c_1(TW^{l_1}(X_t)+1)}{T^{2(k+1)}},\ \ \ t\in [0,T]$$ holds for some constant $c_1>0$ independent of $T$. Thus, for any $\ll>0$,

\beq\label{EE1} \beg{split} \E \e^{ \dd(h)/\ll} &= \E\exp\bigg[\ff 1 \ll\int_0^T\<\dot h_t, \d B_t\>\bigg]\le \bigg(\E\exp\bigg[\ff 2{\ll^2}\int_0^T |\dot h_t|^2\d t\bigg]\bigg)^{1/2} \\
&\le \bigg(\E\exp\bigg[\ff{c_2}{\ll^2}\Big(\ff{\int_0^T W^{2l_1}(X_t)\d t}{T^{4k+2}} +\ff 1 {T^{4k+3}}\Big)\bigg]\bigg)^{1/2}.\end{split}\end{equation}
On the other hand, since $l_1\in [0,1],$ by Lemma \ref{L1} and the Jensen inequality, there exist two constants $c_3,c_4>0$ such that
\beq\label{EE2} \E\exp\bigg[\ff{c_3l_1}T\int_0^TW(X_t)\d t\bigg] \le \e^{c_4l_1 W},\ \ T\in (0,1].\end{equation}
Moreover, since $2l_1<1$, there exists a constant $c_5>0$ such that
$$\ff{c_2  W^{2l_1}}{\ll^2T^{4k+2}}\le\ff{c_3l_1W}{(1+\ll)^2T} + \ff{c_5(1+\ll^{-1})^{4l_1/(1-2l_1)}}{\ll^2T^{(4k+2-2l_1)/(1-2l_1)}},
\ \ \ll,T>0.$$ Combining this with
(\ref{EE1}) and (\ref{EE2}), we conclude that
$$\log \E\e^{\dd(h)/\ll}\le  \ff{cl_1 W}{(1+\ll)^2}+\ff{c(1+\ll^{-1})^{4l_1/(1-2l_1)}}{\ll^2T^{(4k+2-2l_1)/(1-2l_1)}}+\ff{c}{\ll^2 T^{4k+3}}, \ \ T\in (0,1],\ll>0$$ holds for some constant $c>0.$ This completes the proof of (3) by (\ref{BS}) and the Young inequality (see \cite[Lemma 2.4]{ATW09})
\beq\label{LL}|\nn P_T f|= |\E [f(X_T) \dd(h)]|\le \ll \big\{P_Tf\log f-(P_Tf)\log P_T f\big\}+\ll(P_Tf)\log\E\e^{\dd(h)/\ll}.\end{equation}

(4) Again, we only consider $T\in (0,1].$ Let $c_2$ and $C$ be in (\ref{EE1}) and Lemma \ref{L1} respectively.
Then there exists a constant  $c>0$  such that
for any $T\in (0,1]$, $\ll \ge \ff c {T^{2k}}$ implies
$$\ff{c_2}{\ll^2T^{4k+2}}\le \ff 2 {T^2 C \|\si\|^2 \e^{4+2CT}}.$$ Thus, by (\ref{EE1}) and Lemma \ref{L1}, if  $\ll \ge \ff c {T^{2k}}$ then
$$\log \E\e^{\dd(h)/\ll} \le \ff{c_2T^2C\|\si\|^2\e^{4+2CT}}{4\ll^2T^{4k+2}}\log \E\exp\bigg[\ff {2\int_0^T W(X_t) \d t}{T^2 C\|\si\|^2 \e^{4+2CT}}\bigg] +\ff{c_2}{\ll^2T^{4k+3}}
\le \ff{c'(T^2W+1)}{\ll^2T^{4k+3}}$$ holds for some constant $c'>0$ independent of $T$. Combining this with (\ref{LL}) we finish the proof.
\end{proof}

To derive the Harnack inequality of $P_T$ from Theorem \ref{T4.2} (3) and (4), let us  recall a result of \cite{GW}. If there exist a constant
$\ll_0>0$ and
 a positive measurable function $\gg: [\ll_0,\infty)\times \R^{m+d}\to [0,\infty)$ such that
\beq\label{*1} |\nn_v P_Tf|\le \ll \big\{P_Tf\log f-(P_Tf)\log P_Tf\big\} + \gg(\ll,\cdot) P_Tf,\ \ \ll\ge \ll_0\end{equation} holds for some
constant $\ll_0\in (0,\infty]$ and all $f\in \B_b^+(\R^{m+d}),$ then by \cite[Proposition 4.1]{GW},
\beq\label{*2} P_T f(x)\le (P_Tf^p)^{1/p}(x+v)\exp\bigg[ \int_0^1\ff{\gg(\ff{p-1}{1+(p-1)s}, x+sv)}{1+(p-1)s} \d s\bigg]\end{equation}
holds for all $f\in\B_b^+(\R^{m+d})$ and  $p\ge 1+\ll_0.$
Then we have the following consequence of Theorem \ref{T4.2} (3) and (4).

\beg{cor} Let {\bf (H)} and {\bf (A)} hold such that $\nn^{(2)}Z^{(1)}=B_0$ is constant. \beg{enumerate} \item[$(1)$] If $l_1\in [0,1/2)$, then there exists a constant $c>0$ such that
\beg{equation*}\beg{split} &P_Tf(x)\le (P_Tf^p)^{1/p}(x+v)\\
&\times \exp\bigg[ \ff{c|v|^2}{p-1}\Big(\ff{(p-1)l_1\int_0^1 W(x+sv)\d s}{p-1+|v|} +\ff{(1+\ff{p|v|}{p-1})^{4l_1/(1-2l_1)}}
{(T\land 1)^{(4k+2-2l_1)/(1-2l_1)}}+\ff 1 {T^{4k+3}}\Big)\bigg]\end{split}\end{equation*} holds for all $x,v\in\R^{m+d}, T>0, p>1$ and $f\in \B_b^+(\R^{m+d}).$
\item[$(2)$] If $l_1=1$ then there exist two constants $c,c'>0$ such that for any $ T>0,f\in\R^{m+d}$ and $x,v\in \R^{m+d}$,
$$P_Tf(x)\le (P_Tf^p)^{1/p}(x+v)\exp\bigg[\ff{c'|v|^2\big\{1+(T\land 1)^2\int_0^1W(x+sv)\d s\big\}}{(p-1)(T\land 1)^{4k+3}} \bigg]$$
holds for $p\ge 1 +\ff{c|v|}{(T\land 1)^{2k}}.$  \end{enumerate}\end{cor}

\beg{proof} (1) Let $v\in\R^{m+d}$ with $|v|>0$. By Theorem \ref{T4.2}(3), we have
\beg{equation*}\beg{split} |\nn_v P_Tf|\le &\ll |v|\big\{P_Tf\log f-(P_Tf)\log P_T f\big\} \\
&+\ff {c|v|}  \ll \bigg\{\ff{ l_1W}{(1+\ll^{-1})^2}+\ff{(1+\ll^{-1})^{4l_1/(1-2l_1)}}{(T\land 1)^{(4k+2-2l_1)/(1-2l_1)}}
+\ff 1 {(T\land 1)^{4k+3}}\bigg\}P_Tf,\  \ll>0.\end{split}\end{equation*}
Replacing $\ll$ by $\ff \ll {|v|}$, we see that (\ref{*1}) holds for any $\ll_0>0$ and

$$\gg(\ll,\cdot)= \ff {c|v|^2}  \ll \bigg\{\ff{ l_1W}{(1+|v| \ll^{-1})^2}+\ff{(1+|v| \ll^{-1})^{4l_1/(1-2l_1)}}{(T\land 1)^{(4k+2-2l_1)/(1-2l_1)}}+\ff 1 {(T\land 1)^{4k+3}}\bigg\},\ \ \ll>0.$$
Then the desired Harnack inequality follows from (\ref{*2}) since
\beg{equation*}\beg{split} &\int_0^1 \ff{\gg(\ff{p-1}{1+(p-1)s}, x+sv)}{1+(p-1)s} \d s\\
&=\ff{c|v|^2}{p-1} \int_0^1 \bigg\{\ff{l_1W(x+sv)}{1+\ff{|v|(1+(p-1)s)}{p-1}}+\ff{(1+\ff{|v|(1+(p-1)s)}{p-1})^{4l_1/(1-2l_1)}}{(T\land 1)^{(4k+2-2l_1)/(1-2l_1)}}+\ff 1 {(T\land 1)^{4k+3}}\bigg\}\d s\\
&\le \ff{c|v|^2}{p-1}\bigg(\ff{l_1(p-1)\int_0^1W(x+sv)\d s}{p-1+|v|} + \ff{(1+\ff{p |v|}{p-1})^{4l_1/(1-2l_1)}}{(T\land 1)^{(4k+2-2l_1)/(1-2l_1)}}+\ff 1 {(T\land 1)^{4k+3}}\bigg).
\end{split}\end{equation*}

(2) Let $v\in \R^{m+d}$ with $|v|>0.$ By Theorem \ref{T4.2}(4),
$$|\nn_v P_Tf|\le |v| \ll\big\{P_Tf\log f-(P_Tf)\log P_T f\big\}+ \ff{c'|v|((1\land T)^2W+1)}{\ll (T\land 1)^{4k+3}}P_T f$$ holds for
$\ll\ge \ff{c}{(T\land 1)^{2k}}.$
Using $\ff{\ll}{|v|}$ to replace $\ll$, we see that (\ref{*1}) holds for $\ll_0= \ff{c|v|}{(T\land 1)^{2k}}$ and
$$\gg(\ll,\cdot)= \ff{c' |v|^2 ((1\land T)^2W+1)}{\ll (T\land 1)^{4k+3}}.$$Then the proof is completed by (\ref{*2}).
\end{proof}

Finally, according to e.g. \cite[\S4.2]{W11},  the Harnack inequalities presented above imply explicit heat kernel estimates and entropy-cost inequalities for the invariant  probability measure (if exists).

Since  there exist many  non-trivial examples of $A$ and $B_0$ such that {\bf (A)} holds (see \cite{Kal}), it is easy to construct corresponding examples to illustrate Theorem \ref{T4.2}. For instance, for Theorem \ref{T4.2} (3) and (4) only simply consider (\ref{HS}) with $H(x,y)=\<Ax,y\>+W(y)$ such that   $\nn W=B_0$, and for assertion (2) a small perturbation of $W$ is allowed.

\paragraph{Acknowledgement} The authors would like to thank the referee for helpful comments.

\end{document}